# DISCUSSION OF "LEAST ANGLE REGRESSION" BY EFRON ET AL.

By Keith Knight

*University of Toronto*

First, I congratulate the authors for a truly stimulating paper. The paper resolves a number of important questions but, at the same time, raises many others. I would like to focus my comments to two specific points.

**1.** The similarity of Stagewise and LARS fitting to the Lasso suggests that the estimates produced by Stagewise and LARS fitting may minimize an objective function that is similar to the appropriate Lasso objective function. It is not at all (at least to me) obvious how this might work though. I note, though, that the construction of such an objective function may be easier than it seems. For example, in the case of bagging [Breiman (1996)] or subagging [Bühlmann and Yu (2002)], an "implied" objective function can be constructed. Suppose that $\widehat{\theta}_1, \ldots, \widehat{\theta}_m$ are estimates (e.g., computed from subsamples or bootstrap samples) that minimize, respectively, objective functions $Z_1, \ldots, Z_m$ and define

$$\widehat{\theta} = g(\widehat{\theta}_1, \ldots, \widehat{\theta}_m);$$

then $\widehat{\theta}$ minimizes the objective function

$$Z(t) = \inf\{Z_1(t_1) + \cdots + Z_m(t_m) : g(t_1, \ldots, t_m) = t\}.$$

(Thanks to Gib Bassett for pointing this out to me.) A similar construction for stagewise fitting (or LARS in general) could facilitate the analysis of the statistical properties of the estimators obtained via these algorithms.

**2.** When I first started experimenting with the Lasso, I was impressed by its robustness to small changes in its tuning parameter relative to more classical stepwise subset selection methods such as Forward Selection and Backward Elimination. (This is well illustrated by Figure 5; at its best, Forward Selection is comparable to LARS, Stagewise and the Lasso but

---

the performance of Forward Selection is highly dependent on the model size.) Upon reflection, I realized that there was a simple explanation for this robustness. Specifically, the strict convexity in $\boldsymbol{\beta}$ for each $t$ in the Lasso objective function (1.5) together with the continuity (in the appropriate sense) in $t$ of these objective functions implies that the Lasso solutions $\widehat{\boldsymbol{\beta}}(t)$ are continuous in $t$; this continuity breaks down for nonconvex objective functions. Of course, the same can be said of other penalized least squares estimates whose penalty is convex. What seems to make the Lasso special is (i) its ability to produce exact 0 estimates and (ii) the "fact" that its bias seems to be more controllable than it is for other methods (e.g., ridge regression, which naturally overshrinks large effects) in the sense that for a fixed tuning parameter the bias is bounded by a constant that depends on the design but not the true parameter values. At the same time, though, it is perhaps unfair to compare stepwise methods to the Lasso, LARS or Stagewise fitting since the space of models considered by the latter methods seems to be "nicer" than it is for the former and (perhaps more important) since the underlying motivation for using Forward Selection is typically not prediction. For example, bagged Forward Selection might perform as well as the other methods in many situations.

## REFERENCES

Breiman, L. (1996). Bagging predictors. *Machine Learning* **24** 123–140.
Bühlmann, P. and Yu, B. (2002). Analyzing bagging. *Ann. Statist.* **30** 927–961.
  MR1926165

Department of Statistics
University of Toronto
100 St. George St.
Toronto, Ontario M5S 3G3
Canada
e-mail: keith@utstat.toronto.edu